# Smoothing $\ell_1$-penalized estimators for high-dimensional time-course data

**Lukas Meier and Peter Bühlmann**

*Seminar für Statistik, ETH Zürich*
*Leonhardstrasse 27, CH-8092 Zürich, Switzerland*
*e-mail:* `meier@stat.math.ethz.ch`; `buhlmann@stat.math.ethz.ch`

**Abstract:** When a series of (related) linear models has to be estimated it is often appropriate to combine the different data-sets to construct more efficient estimators. We use $\ell_1$-penalized estimators like the Lasso or the Adaptive Lasso which can simultaneously do parameter estimation and model selection. We show that for a time-course of high-dimensional linear models the convergence rates of the Lasso and of the Adaptive Lasso can be improved by combining the different time-points in a suitable way. Moreover, the Adaptive Lasso still enjoys oracle properties and consistent variable selection. The finite sample properties of the proposed methods are illustrated on simulated data and on a real problem of motif finding in DNA sequences.

## 1. Introduction

The Lasso [16] has attracted a lot of attention for prediction and variable selection in linear regression models, including high-dimensional settings where the number of covariates is much larger than sample size [6, 12, 3, 20, 13, 22]. Not only has the idea of $\ell_1$-penalization shown its success in other models [17, 9, 14], but also many extensions of the Lasso in linear regression models have been proposed, among them are the Fused Lasso [18], the Adaptive Lasso [24] and the Relaxed Lasso [11].

Also for multivariate regression, penalization estimators have been shown to be successful [19, 15]. In many problems there is a natural ordering of the response space: our new methodology and theory are exploiting this fact. If we think of time-course data where we observe a response variable over certain time-points, the relationship between "neighbouring" time-points is expected to be stronger than between more distant time-points. Instead of separately estimating a parameter vector for each time-point, it is often a better strategy to combine information across different time-points. By putting an appropriate constraint on the parameter vector, we can control certain characteristics, e.g.

the smoothness over time. As an advantage, we may get a more efficient estimator: By pooling of information, we reduce the variance, typically at the cost of some bias, to achieve a lower mean squared error. For multivariate regression, [2] use the correlation between the responses to construct an estimator with lower mean squared prediction error. In the discussion of [2], the idea of relevance weighted likelihood [7] is mentioned [23]. We use this idea for $\ell_1$-penalized estimators. By using an estimator which also fits well for neighbouring time-points, we can not only get a smoother behaviour of the parameter vector over time, but also profit from more efficiency, both in estimation accuracy and in variable selection.

The rest of this article is organized as follows. In Section 2 we introduce the Smoothed Lasso estimator and show that it asymptotically reduces the bound on the mean squared error compared to the univariate Lasso estimator. In Section 3 we apply the smoothing idea to the Adaptive Lasso and variants thereof and show that it can consistently select the correct model and has a faster convergence rate than the univariate estimator. Simulations follow in Section 4 and a real data analysis for motif search in DNA sequences in Section 5. Section 6 contains the discussion. All proofs can be found in the Appendix.

## 2. Smoothed Lasso

We first start with the definition of the Smoothed Lasso estimator and then study its theoretical properties.

### *2.1. Definition*

Assume that we observe data at $N$ different time-points and that at each time-point $t_r, r = 1, \ldots, N$, we have a linear regression problem of the form

$$y(t_r) = X\beta(t_r) + \varepsilon(t_r),$$

where $X$ is the $n \times p$ design matrix, $y(t_r) \in \mathbb{R}^n$ is the response vector, $\beta(t_r) \in \mathbb{R}^p$ is the parameter vector corresponding to time-point $t_r$ and $\varepsilon(t_r) \in \mathbb{R}^n$ is the corresponding error vector: $\varepsilon(t_r), r = 1, \ldots, N$ are assumed to be i.i.d. random vectors with i.i.d. components having mean zero and finite variance $\sigma^2$. Note that the design matrix $X$ does not depend on $t_r$ in our setup (but it could), and hence we consider a multivariate linear model. As commonly used for penalized estimation, we assume that the columns of $X$ are centered and scaled to have empirical column means 0 and column variances 1.

**Remark 2.1.** *Generalizations of the methodology and theory include that the design matrix $X$ depends on $t_r$, i.e. $X(t_r)$, and that the errors have correlated components $\mathrm{Cov}(\varepsilon(t_r)) = \Sigma$ or arise from a dependent, stationary process with respect to the time-points.*

The idea of the Smoothed Lasso is to use an $\ell_1$-penalty and to suitably combine or smooth the information of the different time-points. It is defined as

$$\hat{\beta}_{\lambda_n,w}(t_r) = \arg\min_\beta \sum_{s=1}^N w(t_s,t_r)\|y(t_s) - X\beta\|_2^2 + \lambda_n \sum_{j=1}^p |\beta_j|, \qquad (2.1)$$

where $w(t_s,t_r)$ are weights satisfying $\sum_{s=1}^N w(t_s,t_r) = 1$. A typical choice is

$$w(t_s,t_r) \propto K\left(\frac{t_s - t_r}{h}\right),$$

where $K(\cdot)$ is a univariate kernel, i.e. $K(x) \geq 0$, $K(x) = K(-x)$, $\int_{-\infty}^\infty K(x)dx = 1$, and $h$ is a bandwidth parameter. Thus, the Smoothed Lasso is $\hat{\beta}(t_r) = \hat{\beta}_{\lambda_n,h}(t_r)$, depending on two tuning parameters.

We can rewrite the weighted optimization problem (2.1) as an ordinary Lasso problem

$$\hat{\beta}_{\lambda_n,h}(t_r) = \arg\min_\beta \|\tilde{y}(t_r) - X\beta\|_2^2 + \lambda_n \sum_{j=1}^p |\beta_j|, \qquad (2.2)$$

where

$$\tilde{y}(t_r) = \sum_{s=1}^N w(t_s,t_r)y(t_s).$$

Hence, any algorithm to solve a standard Lasso problem can be used to calculate the Smoothed Lasso estimator for a given bandwidth $h$.

By forcing the estimate $\hat{\beta}(t_r)$ to fit also well for "neighbouring" time-points, a smooth (non-parametric) trend of $\hat{\beta}(t_r)$ as a function of time is usually inherited (if $p < n$ this is always true because of strict convexity with respect to $\beta$ and continuity with respect to $t_r$ of the criterion in (2.2)).

**Remark 2.2.** *Another approach would be to use a Fused Lasso penalty which also penalizes the absolute value of the differences between neighbouring time-points, i.e. $|\beta_j(t_r) - \beta_j(t_{r-1})|$. Such an approach has two drawbacks: First, we have to model all time-points simultaneously, i.e. fit a model with $Np$ parameters. Moreover, the Fused Lasso problem is more difficult to solve than the Lasso problem. In our approach we fit $N$ Lasso problems with $p$ parameters each. Second, if we want to mimick the behaviour of a bandwidth which is locally adaptive to the underlying true parameter function $\beta(t)$, we have to introduce a lot of tuning parameters for the Fused Lasso and search over a high-dimensional grid when doing cross-validation.*

**Remark 2.3.** *We do not assume that the active set (the set of predictors with nonzero coefficients) stays the same over all time-points. Our methodology allows for the fact that some predictors enter or leave the active set along the time-course.*

In the next Section we first consider the special case of an orthogonal design and indeed, we show that the mean squared error is decreased asymptotically.

## 2.2. Orthogonal Case

We consider the situation where the number of parameters equals the number of observations and the design matrix is orthogonal, i.e. $X^T X = nI_n$, and where the errors $\varepsilon(t_r)$ are Gaussian. We focus on a single time-point of interest and therefore omit the time-index for notational simplicity. In [5, Theorem 1] it is shown that the (univariate) soft-threshold estimator $\hat{\beta}_{ST}$ (with threshold parameter $\sigma\sqrt{2\log(n)/n}$), which is equivalent to the Lasso in the orthogonal case (with penalty parameter $\lambda = 2\sigma\sqrt{2\log(n)n}$), satisfies

$$\mathbb{E}[\|\hat{\beta}_{ST} - \beta\|_2^2] \leq (2\log(n) + 1)\left\{\frac{\sigma^2}{n} + \sum_{i=1}^n \min\left(\beta_i^2, \frac{\sigma^2}{n}\right)\right\} \quad (2.3)$$

for all $\beta \in \mathbb{R}^n$ and that this bound is asymptotically sharp in a minimax-sense [5, Theorem 3]. If the non-zero $\beta_i$'s are not of too low order (i.e. $|\beta_i| \gg n^{-1/2}$), the order of this bound is $\log(n)|\mathcal{A}_n|/n$, where $\mathcal{A}_n = \{i \,|\, \beta_i \neq 0\}$ denotes the active set of the time-point of interest. Even though we restrict ourselves to a class of parameter vectors which stay out of the $n^{-1/2}$-range, the order of the bound in (2.3) is sharp because the maximal risk is attained for an element of this class (see the proof of Theorem 3 in [5]).

The order $\log(n)|\mathcal{A}_n|/n$ can be decreased by the smoothed estimator as shown in Proposition 2.1.

**Proposition 2.1.** *Assume Gaussian errors $\varepsilon(t_r)$ and the regularity conditions (RC 1) – (RC 4) described in Appendix A.1. For $h = h_n \asymp \log(n)^{1/5} n^{-1/5} N^{-1/5}$ and $N = N_n$ such that $Nh \to \infty$ $(n \to \infty)$ the risk of the Smoothed Lasso from (2.2) asymptotically satisfies: for $\lambda_n = 2\sigma(Nh)^{-1/2}\sqrt{2\log(n)n}$*

$$\begin{aligned}\mathbb{E}[\|\hat{\beta}_{\lambda_n,h_n} - \beta\|_2^2] &\leq C\log(n)|\mathcal{A}_n|/(nNh) \\ &\asymp \log(n)^{4/5}|\mathcal{A}_n|/(n^{4/5}N^{4/5}),\ n \to \infty.\end{aligned}$$

*for all $\beta \in \mathbb{R}^n$ and some constant $C$.*

A proof is given in Appendix A.2.

For a faster convergence rate than $\log(n)|\mathcal{A}_n|/n$ (for the unsmoothed Lasso) we require $Nh$ to converge to infinity which implies that

$$N = N_n \gg \left(\frac{n}{\log(n)}\right)^{1/4},$$

i.e. $N$ can be of much lower order than $n$ for achieving a faster convergence rate for the minimax bound.

## 2.3. General Case

Let us now consider the general case, i.e. we do not restrict ourselves to an orthogonal design matrix. In particular, we allow for high-dimensional situations where $p = p_n \gg n$ is increasing very fast as $n \to \infty$.

Using the results in [13] for a fixed design, the univariate Lasso estimator satisfies under regularity conditions on the design matrix

$$\|\hat{\beta}_{Lasso,\lambda_n} - \beta\|_2^2 \leq O_P\left(\sigma^2 \frac{m_{\lambda_n}}{n} \log(p_n)\right) + O\left(\frac{|\mathcal{A}_n|}{m_{\lambda_n}}\right), \; n \to \infty,$$

where $m_{\lambda_n} = Cn^2/\lambda_n^2$ for some constant $C$. In a certain sense, this bound is tight, see [13, Remark 1]. We can choose

$$\lambda_n \asymp \sigma^{1/2} n^{3/4} \log(p_n)^{1/4} |\mathcal{A}_n|^{-1/4}$$

and arrive at the optimal rate

$$\|\hat{\beta}_{Lasso,\lambda_n} - \beta\|_2^2 \leq O_P\left(\sigma n^{-1/2} \log(p_n)^{1/2} |\mathcal{A}_n|^{1/2}\right), \; n \to \infty. \qquad (2.4)$$

For proving such a result, assumptions on the design matrix are crucial: Various authors use different conditions, cf. [13, 3, 20, 22]. We refer the reader to [13] for a detailed description of the regularity conditions for (2.4).

**Proposition 2.2.** *Assume that the univariate Lasso satisfies (2.4) and denote the bound on the right hand side of (2.4) by $a_n = \sigma n^{-1/2} \log(p_n)^{1/2} |\mathcal{A}_n|^{1/2}$. Furthermore, assume the regularity conditions (RC1) – (RC4) described in Appendix A.1. Then, if $N = N_n \gg |\mathcal{A}_n|^{1/4} a_n^{-1/4}$ and for some suitable $\lambda_n$ and $h = h_n$:*

$$\|\hat{\beta}_{\lambda_n, h_n} - \beta\|_2^2 = o_P(a_n),$$

*i.e. the Smoothed Lasso has a faster convergence rate than the (tight) bound in (2.4) for the Lasso.*

A proof is given in Appendix A.2. Suitable choices for $\lambda_n$ and $h_n$ in Proposition 2.2 are

$$h_n \asymp N^{-1/9} |\mathcal{A}_n|^{-2/9} a_n^{2/9}$$

and

$$\lambda_n \asymp \sigma^{1/2} (Nh)^{-1/4} n^{3/4} \log(p_n)^{1/4} |\mathcal{A}_n|^{-1/4}.$$

Using the notation that is introduced at the beginning of the proof of Proposition 2.1, one can derive other asymptotic properties by linking known results for the Lasso [6, 3, 20] with the smoothed model

$$\tilde{y} = X\tilde{\beta} + \tilde{\varepsilon}$$

and an analysis of the bias term $\|\tilde{\beta} - \beta\|_q$ for $q \in \{1, 2\}$ as in (A.3).

Up to now we only considered the estimation error for $\beta$ and no variable selection properties. The smoothing reduces the variance and thus it can be expected that the Smoothed Lasso selects more (noise) variables than its univariate counterpart. Empirical evidence of this property is given in Section 4. This problem can be overcome by a second stage which removes many of the coefficients whose estimates are close to zero. In fact, already the case with a univariate response often requires such a second stage for consistent variable selection [24]. We will treat a special case in the next Section.

## 3. Smoothed Adaptive Lasso

The Adaptive Lasso [24] weights the penalty for the different coefficients using an initial estimator $\hat{\beta}_{init}$, i.e.

$$\hat{\beta}^{(\hat{\beta}_{init})}_{\lambda_n} = \arg\min_{\beta} \|y - X\beta\|_2^2 + \lambda_n \sum_{j=1}^{p} \hat{\tau}_j |\beta_j|,$$

where $\hat{\tau}_j = 1/|\hat{\beta}_{init,j}|^\gamma$ for some $\gamma > 0$ are weights based on the initial estimator $\hat{\beta}_{init}$. For simplicity we will restrict ourselves to $\gamma = 1$. In [24], the ordinary least squares (OLS) estimator is used for $\hat{\beta}_{init}$: here, we will mainly use the Lasso and versions thereof. Through a re-scaling of the columns of the design matrix, the Adaptive Lasso estimator can be formulated as an ordinary Lasso problem, see [24].

We can also apply the smoothing technique of Section 2 to the Adaptive Lasso. In the smoothed case we again replace the residual sum of squares in the objective function with its smoothed counterpart in (2.2), i.e.

$$\hat{\beta}^{(\hat{\beta}_{init})}_{\lambda_n,h}(t_r) = \arg\min_{\beta} \|\tilde{y}(t_r) - X\beta\|_2^2 + \lambda_n \sum_{j=1}^{p} \hat{\tau}_j |\beta_j|. \qquad (3.1)$$

In [24], an asymptotic oracle result for the Adaptive Lasso is given for fixed dimension $p$. We show that the Smoothed Adaptive Lasso has a faster convergence rate. Again, as we focus on a single time-point, we omit the time-index for notational simplicity.

We will consider the situation where the number of variables $p$ is kept fixed as $n \to \infty$. As before, let $\mathcal{A}$ be the active set of the true parameter vector at the current time-point and $\hat{\mathcal{A}}_n$ be its empirical counterpart.

**Theorem 3.1.** *Assume a fixed design with $\lim_{n\to\infty} \frac{1}{n} X^T X = C$ for some positive definite matrix $C$. If $\hat{\beta}_{init} - \beta = O_P(a_n^{-1})$ for some sequence $a_n \to \infty$, $\lambda_n\sqrt{N_n h_n}/\sqrt{n} \to 0$, $\lambda_n a_n \sqrt{N_n h_n}/\sqrt{n} \to \infty$, $h_n = o(n^{-1/5} N_n^{-1/5})$ and $N_n h_n \to \infty$ $(n \to \infty)$, then the Smoothed Adaptive Lasso in (3.1) satisfies under the regularity conditions (RC1) – (RC4) described in Appendix A.1:*

$$\lim_{n\to\infty} \mathbb{P}(\hat{\mathcal{A}}_n = \mathcal{A}) = 1$$

*and*

$$\sqrt{nN_n h_n}(\hat{\beta}^{(\hat{\beta}_{init})}_{\lambda_n,h_n,\mathcal{A}} - \beta_{\mathcal{A}}) \xrightarrow{d} N(0, \sigma_*^2 (C_{\mathcal{A}\mathcal{A}})^{-1}), \ n \to \infty,$$

*where $\sigma_*^2 = \sigma^2 \int_{-\infty}^{\infty} K^2(x)dx$, $\hat{\beta}_{\mathcal{A}}$, $\beta_{\mathcal{A}}$ are the sub-vectors of $\hat{\beta}$, $\beta$ and $C_{\mathcal{A}\mathcal{A}}$ is the sub-matrix of $C$ corresponding to the active set $\mathcal{A}$.*

A proof is given in Appendix A.2. Thus, if the initial estimator is consistent, we can find a sequence $\lambda_n$ such that the Smoothed Adaptive Lasso has the property of consistent model selection and asymptotic normality on the active set $\mathcal{A}$. Furthermore, if $N = N_n \gg n^{1/4}$, we can choose $h = h_n = o(n^{-1/5}N_n^{-1/5})$ such that $Nh \to \infty$. Thus, as already pointed out in Section 2.2, a relatively small value of $N = N_n$ is sufficient for achieving an improved convergence rate.

**Remark 3.1.** *The optimal convergence rate $n^{-2/5}N^{-2/5}$ in Theorem 3.1 can be achieved using $h_n \asymp n^{-1/5}N^{-1/5}$. Then, the limiting normal distribution becomes $\mathcal{N}(B_\mathcal{A}, \sigma_*^2(C_{\mathcal{A}\mathcal{A}})^{-1})$ for some vector $B_\mathcal{A}$ with $0 \le |B_{\mathcal{A},j}| < \infty$ for all $j$. This is the same distribution as when using local least squares (with kernel $K$). Hence, the Smoothed Adaptive Lasso has an oracle property saying that it is asymptotically as good as local least squares with the true underlying active set $\mathcal{A}$ known beforehand.*

### 3.1. Choice of initial estimator

The choice of the initial estimator will influence the final estimator. In particular, the sparsity of the final estimator can be maximized by making an appropriate choice, as discussed below.

We will first focus on univariate estimators, i.e. on estimators which only use the data of the current time-point. In view of Theorem 3.1, the basic assumption for the initial estimator is consistency. The ordinary least squares (OLS) method is a possible choice for low-dimensional problems with fixed dimension $p$ as it is $\sqrt{n}$-consistent. The Lasso is consistent in an $\ell_2$-sense, even in the high-dimensional setting, see Section 2.3. Finally, the Adaptive Lasso is $\sqrt{n}$-consistent for fixed $p$ [24] and consistent under suitable regularity conditions for $p \gg n$ [8]. For high-dimensional problems the OLS estimator is not appropriate because it is unstable or even not defined in a $p > n$ situation. The Lasso or Adaptive Lasso are more appropriate choices.

If the initial estimator is doing variable selection, i.e. some of the coefficients $\hat{\beta}_{init,j} = 0$, the smoothed estimator is at least as sparse as the initial estimator: a zero-coefficient in the initial estimator, i.e. $\hat{\beta}_{init,j} = 0$, results in an infinite penalty for that component, i.e. $\hat{\tau}_j = \infty$, forcing the smoothed estimate to be zero, i.e. $\hat{\beta}_j(t_r) = 0$. This reduces the computational complexity for the smoothing stage since some or even many predictors can be excluded from the model.

For the case that the initial estimator has a tuning parameter, as with the Lasso and the Adaptive Lasso, one would in practice tune it to be prediction optimal. For the Lasso, this produces too large models, i.e. many noise variables are included in the selected model [12]. However, noise variables tend to have small coefficients and will therefore be heavily penalized in the second smoothing step of the Smoothed Adaptive Lasso.

It is of course also possible to use a smoothed estimator as initial estimator, e.g. the Smoothed Lasso. In terms of the number of selected variables, as we will see in Section 4, this is often worse than directly using the univariate counterpart. Due to the reduced variance, the smoothed initial estimator tends to select too many variables and not all of them will be eliminated in the second stage of the Smoothed Adaptive Lasso.

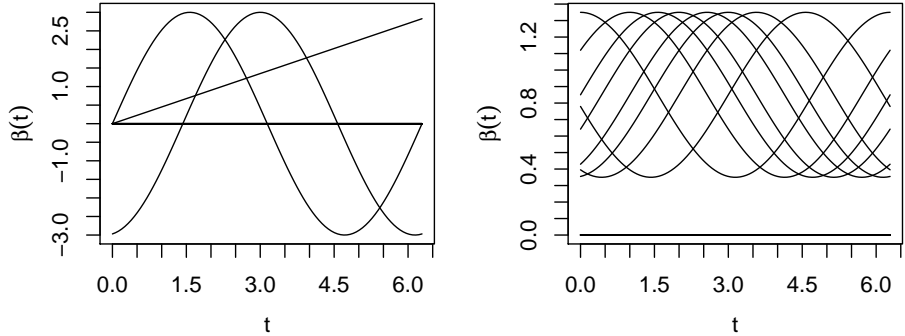

Fig 1. *Parameter functions for model 1 (left) and model 2 (right).*

In view of some empirical results in Section 4, we advocate the following: the initial estimator for the Smoothed Adaptive Lasso is the univariate Adaptive Lasso; the latter itself uses the univariate Lasso as initial estimator. This amounts to be a three-stage procedure where all of the estimations are tuned to be prediction optimal using e.g. some cross-validation scheme. There is substantial agreement by now that two or more stages are needed to achieve good regularization properties in high-dimensional settings [24, 11, 25, 13, 21, 10]. As a novelty here, our third stage involves an additional smoothing operation.

## 4. Simulations

In this Section we want to evaluate the finite sample properties of the proposed estimators.

### 4.1. Design

We consider the following models, similar to [24]:

Model 1: Some large effects

$$\beta(t) = (0.45t,\, 3\sin(t),\, 3\cos(t-3),\, 0,\ldots, 0)$$

Model 2: Many small effects

$$\begin{aligned}\beta(t) &= (0.85 + 0.5\sin(t),\, 0.85 + 0.5\cos(t),\, 0.85 + 0.5\sin(t-1),\\ & \quad 0.85 + 0.5\cos(t-1),\ldots,\, 0.85 + 0.5\sin(t-3),\\ & \quad 0.85 + 0.5\cos(t-3),\, 0,\ldots, 0)\end{aligned}$$

Figure 1 illustrates the two parameter vectors as a function of time $t$. We use an equidistant grid on the interval $[0, 2\pi]$, i.e. $t_r = (r-1)\frac{2\pi}{N-1}, r \in \{1, \ldots, N\}$, where $N = 18$. The design matrix $X$ is simulated from a multivariate normal distribution with mean zero and covariance matrix $\Sigma_{i,j} = 0.5^{|i-j|}$. The standard deviation of the error term is chosen from $\sigma \in \{2, 4\}$ which corresponds to a signal-to-noise ratio (averaged over $N$) of approximately $\{2.7, 0.7\}$ and $\{3.8, 0.9\}$ for model 1 and model 2, respectively. We use both a "classical" setup with $n = 50$, $p = 8$ and a high-dimensional setup with $n = 100$, $p = 1000$.

The best combination of bandwidth $h$ and penalization parameter $\lambda$ is being searched on a two-dimensional grid using an independent validation set of half the size of the training set. This is done independently for each time-point which means that we allow for a locally varying bandwidth. The density of the standard normal distribution is used as kernel function $K(\cdot)$ for the weight function $w(\cdot, \cdot)$, see Section 2.1.

For the (Smoothed) Adaptive Lasso with (Smoothed) Lasso as initial estimator, we first determine the optimal penalization parameter for the initial estimator and keep it fixed when searching for the optimal penalization parameter and bandwidth for the final estimator.

All estimators which we compare are listed in Table 1.

TABLE 1
*Different estimators*

| Univariate Estimators | Smoothed Estimators |
| --- | --- |
| 1. Lasso | 4. Lasso |
| 2. Adapt. Lasso with OLS | 5. Adapt. Lasso with smoothed OLS |
| 3. Adapt. Lasso with univ. Lasso | 6. Adaptive Lasso with 4. |
| | 7. Adaptive Lasso with 3. |

## 4.2. Performance Measures

To measure the goodness of fit and the ability to pick the model of the correct size we define the following performance measures.

For the mean squared error we use

$$MSE_\beta = \frac{1}{N} \sum_{r=1}^{N} \|\hat{\beta}(t_r) - \beta(t_r)\|_2^2.$$

Moreover, we also report the mean squared prediction error for the regression function $x^T \beta(t_r)$

$$MSE_P = \frac{1}{N} \left( \sum_{r=1}^{N} (\hat{\beta}_0(t_r) - \beta_0(t_r))^2 + (\hat{\beta}(t_r) - \beta(t_r))^T \Sigma (\hat{\beta}(t_r) - \beta(t_r)) \right),$$

where $\hat{\beta}_0(t_r)$ is the intercept term (and $\beta_0(t_r) = 0$ for our simulations) and $\Sigma$ is the covariance matrix of the covariates.

For the number of variables we define the mean model size

$$MSize = \frac{1}{N} \sum_{r=1}^{N} |\hat{\mathcal{A}}(t_r)|$$

and the mean number of false positives

$$FP = \frac{1}{N} \sum_{r=1}^{N} \sum_{j=1}^{p} 1_{[\hat{\beta}_j(t_r) \neq 0]} 1_{[\beta_j(t_r) = 0]}.$$

In applied sciences where (possibly expensive) experiments are conducted to verify the selected variables (e.g. in biology), the number of false positives is a crucial quantity one wants to minimize in order to keep the costs low.

### 4.3. Results

The results can be found in Table 2. For the high-dimensional setting we did not consider OLS initial estimators. Several conclusions can be made. Let us first focus on the Lasso estimator. In all simulation settings, smoothing improves the $MSE_\beta$ score substantially. The downside for the Smoothed Lasso estimator is that due to the decreased variance, more noise variables tend to enter the model which results in larger selected models (with more false positives) than for the univariate Lasso estimator. However, in practice one would assign a variable importance score to each coefficient and therefore concentrate first on those with the largest contributions, whereas many of the false positives have small importance scores only.

Also for the Adaptive Lasso, the $MSE_\beta$ scores get decreased by smoothing in all simulation settings. Using a smoothed initial estimator leads to too large models. Take for example Adaptive Lasso with Smoothed Lasso as initial estimator, i.e. proposal 6 in Table 1. As we have described above, the Smoothed Lasso tends to select a too large initial model. Although the Adaptive Lasso can eliminate most noise variables in the second stage due to their large weights from small coefficients of the initial estimator, the resulting models still get a bit too large. However, the estimator is very competitive with respect to prediction performance.

Using a univariate initial estimator, i.e. our proposal 7 in Table 1, to get more reasonably sized models seems to be a good compromise. It does not only produce the sparsest models but is often also competitive with respect to $MSE_\beta$ and $MSE_P$.

### 5. Real data: Motif Finding in DNA Sequences

We apply the smoothing methodology to a problem of motif regression [4]. A motif (typically a 5–15 letter word consisting of letters A, C, G, and T) is a candidate of a binding site of some functional element, e.g. a transcription factor (a protein which regulates gene expression). In [1] a collection of various gene expression time-course experiments and a set of candidate motifs for yeast is provided. Gene expression values for a total of 2587 genes are available and $p = 666$ motif candidates are used to build the motif scores for each gene. These measure how well the motifs are represented in the upstream regions of the genes. We focus on a time-course experiment spanning $N = 12$ different time-points. In summary we have 2587 observations of a 666 dimensional predictor (the motif scores) and a one-dimensional response (the gene expression value) at each of the 12 time-points. Thus, each row of the design matrix $X$ corresponds to a gene and each column to a motif score. The element $x_{i,j}$ measures how well the $j$th motif score is represented in the upstream region of the $i$th gene.

TABLE 2. *Mean values of the different performance measures based on 100 simulation runs for $n = 50$, $p = 8$ (top) and $n = 100$, $p = 1000$ (bottom). Standard deviations are given in parentheses. The low-dimensional case of model 2 can't have false positives because all variables are active.*

|  | $MSE_\beta$ | $MSE_P$ | $MSize$ | $FP$ | $MSE_\beta$ | $MSE_P$ | $MSize$ | $FP$ |
|---|---|---|---|---|---|---|---|---|
| Model 1 | | $\sigma = 2$ | | | | $\sigma = 4$ | | |
| 1. | 0.83 (0.17) | 0.65 (0.10) | 5.67 (0.50) | 2.96 (0.47) | 3.14 (0.67) | 2.50 (0.41) | 5.10 (0.57) | 2.69 (0.49) |
| 2. | 0.71 (0.17) | 0.57 (0.10) | 4.39 (0.50) | 1.82 (0.46) | 3.03 (0.71) | 2.43 (0.42) | 4.13 (0.51) | 1.89 (0.44) |
| 3. | 0.69 (0.15) | 0.56 (0.09) | 4.30 (0.48) | 1.72 (0.44) | 3.05 (0.72) | 2.43 (0.44) | 4.02 (0.55) | 1.79 (0.46) |
| 4. | 0.54 (0.14) | 0.43 (0.09) | 6.26 (0.49) | 3.48 (0.48) | 1.93 (0.50) | 1.51 (0.32) | 6.09 (0.61) | 3.42 (0.56) |
| 5. | 0.46 (0.13) | 0.38 (0.08) | 4.85 (0.50) | 2.15 (0.48) | 1.75 (0.46) | 1.39 (0.30) | 4.89 (0.59) | 2.34 (0.54) |
| 6. | 0.47 (0.13) | 0.38 (0.09) | 4.80 (0.51) | 2.11 (0.49) | 1.80 (0.48) | 1.43 (0.33) | 4.83 (0.61) | 2.29 (0.54) |
| 7. | 0.51 (0.14) | 0.41 (0.09) | 4.06 (0.47) | 1.50 (0.44) | 2.29 (0.68) | 1.73 (0.45) | 3.78 (0.52) | 1.58 (0.43) |
| Model 2 | | $\sigma = 2$ | | | | $\sigma = 4$ | | |
| 1. | 1.10 (0.19) | 0.86 (0.13) | 7.60 (0.17) | – | 3.22 (0.42) | 3.06 (0.62) | 6.52 (0.40) | – |
| 2. | 1.33 (0.23) | 0.97 (0.15) | 7.35 (0.25) | – | 4.10 (0.56) | 3.62 (0.71) | 5.90 (0.52) | – |
| 3. | 1.30 (0.22) | 0.95 (0.15) | 7.26 (0.24) | – | 3.99 (0.53) | 3.54 (0.70) | 5.79 (0.47) | – |
| 4. | 0.51 (0.12) | 0.42 (0.08) | 7.92 (0.08) | – | 1.42 (0.34) | 1.29 (0.33) | 7.62 (0.22) | – |
| 5. | 0.60 (0.14) | 0.47 (0.10) | 7.64 (0.21) | – | 1.72 (0.43) | 1.46 (0.37) | 7.11 (0.35) | – |
| 6. | 0.64 (0.14) | 0.49 (0.10) | 7.60 (0.21) | – | 1.84 (0.45) | 1.53 (0.40) | 6.99 (0.39) | – |
| 7. | 0.84 (0.20) | 0.61 (0.15) | 7.14 (0.28) | – | 2.94 (0.57) | 2.54 (0.79) | 5.62 (0.48) | – |
| Model 1 | | $\sigma = 2$ | | | | $\sigma = 4$ | | |
| 1. | 2.16 (0.36) | 1.48 (0.23) | 27.35 (5.13) | 24.92 (5.11) | 5.67 (0.70) | 4.22 (0.61) | 19.83 (4.86) | 18.05 (4.80) |
| 3. | 1.07 (0.30) | 0.76 (0.18) | 7.11 (1.93) | 4.79 (1.89) | 4.73 (0.84) | 3.36 (0.59) | 6.07 (2.08) | 4.51 (2.02) |
| 4. | 1.05 (0.18) | 0.78 (0.12) | 41.51 (6.62) | 38.87 (6.61) | 3.12 (0.50) | 2.41 (0.40) | 42.94 (7.01) | 40.54 (6.99) |
| 6. | 0.48 (0.11) | 0.40 (0.09) | 12.74 (3.44) | 10.14 (3.42) | 2.25 (0.50) | 1.81 (0.33) | 17.66 (4.34) | 15.36 (4.34) |
| 7. | 0.74 (0.29) | 0.49 (0.14) | 5.10 (1.32) | 2.79 (1.28) | 3.80 (0.90) | 2.47 (0.57) | 4.74 (1.53) | 3.19 (1.47) |
| Model 2 | | $\sigma = 2$ | | | | $\sigma = 4$ | | |
| 1. | 1.37 (0.23) | 1.69 (0.34) | 32.90 (6.13) | 25.71 (6.17) | 3.93 (0.56) | 5.42 (0.95) | 28.24 (5.06) | 22.84 (5.20) |
| 3. | 1.52 (0.22) | 1.22 (0.23) | 12.13 (2.16) | 5.82 (2.14) | 5.43 (0.57) | 5.35 (0.81) | 11.70 (2.76) | 7.54 (2.70) |
| 4. | 0.50 (0.08) | 0.51 (0.09) | 46.56 (8.63) | 38.60 (8.62) | 1.41 (0.26) | 1.49 (0.31) | 51.24 (8.52) | 43.49 (8.51) |
| 6. | 0.51 (0.09) | 0.44 (0.08) | 23.41 (5.70) | 15.59 (5.70) | 1.72 (0.31) | 1.55 (0.31) | 28.61 (6.66) | 21.22 (6.62) |
| 7. | 1.00 (0.19) | 0.67 (0.14) | 9.74 (1.40) | 3.52 (1.35) | 3.90 (0.50) | 3.23 (0.60) | 9.14 (1.99) | 5.07 (1.92) |

To illustrate the smoothing methods and the effect of different sizes for the training set, we use random subsets of different sizes as training set. An independent validation set is used to determine the prediction optimal tuning parameters. The size of the validation set is half the size of the training set. The remaining data is used as test-set.

The results for a training set of size 1300 is given in Table 3. In terms of prediction error, there is not much gain when smoothing the estimators for this data-set, especially for the Adaptive Lasso. A reason for this may be the large variance of the error term. Note that for a new test observation $(x_{test}, y_{test})$ we have

$$\mathbb{E}_{x_{test}, y_{test}}[(\hat{y} - y_{test})^2] = \mathbb{E}[(x^T\hat{\beta} - x^T\beta)^2] + \sigma^2.$$

The error variance $\sigma^2$ is likely to be the dominating quantity since motif regression is known to be very noisy. In terms of variable selection, the smooth-ing step decreases the model size for the Adaptive Lasso estimator and is potentially reducing the number of false positives: In particular for time-points $t_r = 1, 3, 5, 7, 8, 9, 10$ the Smoothed Adaptive Lasso yields much sparser model fits. For the Lasso estimator, the smoothing has a tendency to increase the number of selected variables resulting in rather large models. This coincides with our findings in Section 4. If we decrease the training sample size to 200 we see some small improvement with respect to the mean squared prediction error (not shown).

TABLE 3
*Mean squared prediction error (top) and number of selected variables (bottom) for the training set of size 1300. Rows: 4 different methods, as described in Table 1. Columns: 12 different time-points.*

|    | 1    | 2    | 3    | 4    | 5    | 6    | 7    | 8    | 9    | 10   | 11   | 12   |
|----|------|------|------|------|------|------|------|------|------|------|------|------|
| 1. | 0.07 | 0.13 | 0.44 | 2.71 | 1.52 | 1.66 | 1.98 | 2.55 | 3.17 | 2.77 | 3.02 | 2.92 |
| 3. | 0.07 | 0.13 | 0.45 | 2.82 | 1.59 | 1.75 | 2.05 | 2.64 | 3.26 | 2.81 | 3.09 | 3.04 |
| 4. | 0.07 | 0.13 | 0.44 | 2.71 | 1.52 | 1.66 | 1.99 | 2.52 | 3.17 | 2.78 | 3.02 | 2.91 |
| 7. | 0.07 | 0.13 | 0.45 | 2.86 | 1.59 | 1.76 | 2.06 | 2.65 | 3.27 | 2.85 | 3.10 | 3.06 |
| 1. | 35   | 28   | 74   | 155  | 133  | 157  | 154  | 141  | 124  | 106  | 123  | 109  |
| 3. | 10   | 6    | 29   | 80   | 85   | 35   | 47   | 31   | 40   | 30   | 42   | 49   |
| 4. | 35   | 29   | 70   | 155  | 123  | 157  | 177  | 178  | 124  | 98   | 126  | 112  |
| 7. | 4    | 5    | 21   | 71   | 46   | 29   | 34   | 22   | 29   | 22   | 35   | 43   |

## 6. Discussion

We propose smoothing techniques for $\ell_1$-penalized (Lasso-type) estimators for a time-course of high-dimensional linear models. We show theoretically that for the Lasso and the Adaptive Lasso, better estimates in terms of the mean squared error can be obtained by combining the responses of different time-points in a suitably weighted way. Empirically, the Smoothed Adaptive Lasso estimator yields the sparsest models with competitive mean squared error performance when using the univariate Adaptive Lasso as initial estimator. The Smoothed Lasso estimator has very good performance with respect to the mean squared error but selects too many noise variables in general. An additional thresholding stage would be necessary if the primary interest is in variable selection.

The smoothing methodology can also be applied to generalized linear models (GLM). The main difference is that we can't rewrite the smoothed estimator as an ordinary lasso problem as in (2.2). This implies that the computational burden increases: In the worst case (depending on the support of the kernel and the bandwidth $h$), by stacking the response variables and design matrices of the different time-points, the total sample size is $Nn$, which can be substantially larger than $n$ in (2.2), while the dimensionality is still $p$.

Our methodology applies to more general problems than time-course settings. For example, we can directly treat the situation of different (heterogeneous) data-sets $(y(t), X)$, $t = 1, \ldots, N$ (or $(y(t), X(t))$, $t = 1, \ldots, N$) with $n(t) \times 1$ response vectors and $n(t) \times p$ design matrices, where $t$ is the index for the various data-sets. All we need is a suitable pseudo-distance $d(t, s)$ among the different data-sets indexed by $t$ and $s$. The weights in (2.2) are then of the form

$$w(t, s) \propto K\left(\frac{d(t, s)}{h}\right).$$

The pseudo-distance $d(\cdot, \cdot)$ could be learned from the data, e.g. based on some pseudo-metrics for clustering different data-sets.

Whether the multivariate view over different time-points (or different data-sets) pays off for a particular problem is not clear a-priori. However, our methodology encompasses the univariate Lasso methods, by choosing the bandwidth $h = 0$. Hence, using some cross-validation scheme enables to find out whether pooling information over different time-points (or data-sets) is worthwhile and if so, the Smoothed (Adaptive) Lasso from the multivariate approach renders more accurate estimates.

## Acknowledgments

We thank the associate editor and two anonymous referees for constructive comments and Yuan Yuan for providing us with the motif scores data-set.

**Appendix A: Proofs**

*A.1. Regularity Conditions*

We denote by $\beta(t) \in \mathbb{R}^p$ the true underlying parameter vector as a function of time $t$.

(RC1) **Curvature of underlying function**
$\beta_j(t)$ is twice continuously differentiable with $\sup_j \left|\beta_j''(t)\right| \leq C < \infty$ for all $t$ and some constant $C$.

(RC2) **Equidistant grid**
For the asymptotic implications we assume that we have an equidistant grid around the time-point of interest $t_r$ of the form
$$t_s = t_r + \frac{s}{N},$$
where $s = -\lfloor N/2 \rfloor, \ldots, \lfloor N/2 \rfloor$. Note that we enumerate using negative values of $s$ as well.

(RC3) **Sampling Points**
For the time-point of interest $t_r$ we assume that if $\beta_j(t_r) = 0$ it follows that there is an open neighbourhood $U_j \ni t_r$, such that $\beta_j(u) = 0 \ \forall\ u \in U_j$. Moreover, we require $\inf_j \operatorname{diam}(U_j) > \delta$ for some $\delta > 0$. I.e. for the time-point of interest no variable enters or leaves the active set.

(RC4) **Compact kernel**
The kernel function $K(\cdot)$ is assumed to have compact support on $[-1, 1]$.

*A.2. Proofs*

**Proof of Proposition 2.1**
As we focus on a single time-point $t_r$ we omit the time index for notational simplicity, whenever possible. For the smoothed response $\tilde{y}$ we have the model
$$\tilde{y} = X\tilde{\beta} + \tilde{\varepsilon}_N,$$
at the time-point of interest, where
$$\tilde{\beta} = \sum_{s=-\lfloor Nh \rfloor}^{\lfloor Nh \rfloor} w(t_s, t_r) \beta(t_s)$$
and
$$\tilde{\varepsilon}_N = \sum_{s=-\lfloor Nh \rfloor}^{\lfloor Nh \rfloor} w(t_s, t_r) \varepsilon(t_s).$$
Note that $(\tilde{\varepsilon}_N)_1, \ldots, (\tilde{\varepsilon}_N)_n$ are i.i.d. with mean zero and variance given below.

We can now use the decomposition

$$\|\hat{\tilde{\beta}} - \beta\|_2^2 \leq 2\|\hat{\tilde{\beta}} - \tilde{\beta}\|_2^2 + 2\|\tilde{\beta} - \beta\|_2^2. \tag{A.1}$$

The first term is "classical". We can use the theory of [5] with respect to an error term with reduced variance. For the asymptotic variance we have

$$\begin{aligned}
\operatorname{Var}((\tilde{\varepsilon}_N)_i) &= \operatorname{Var}\left(\sum_{s=-\lfloor Nh \rfloor}^{\lfloor Nh \rfloor} w(t_s, t_r)\,\varepsilon(t_s)_i\right) \\
&= \sigma^2 \sum_{s=-\lfloor Nh \rfloor}^{\lfloor Nh \rfloor} w^2\left(t_r + \frac{s}{N}, t_r\right) \\
&= \sigma^2 \frac{\left(\frac{1}{Nh}\right)^2 \sum_{s=-\lfloor Nh \rfloor}^{\lfloor Nh \rfloor} K^2\left(\frac{s/N}{h}\right)}{\left\{\frac{1}{Nh} \sum_{s=-\lfloor Nh \rfloor}^{\lfloor Nh \rfloor} K\left(\frac{s/N}{h}\right)\right\}^2}.
\end{aligned}$$

Using a Riemann sum approximation, we arrive at

$$\operatorname{Var}((\tilde{\varepsilon}_N)_i) \sim \frac{\sigma^2}{Nh}\int K^2(x)dx, \tag{A.2}$$

i.e. the error variance is of order $1/(Nh)$.

Let us now consider the bias term. If $\beta_i(t_r) = 0$ it follows with the compactness assumption of the kernel and (RC3) that for $h = h_n$ small enough $\tilde{\beta}_i(t_r) = 0$. If $\beta_i(t_r) \neq 0$ we have

$$\begin{aligned}
\tilde{\beta}_i(t_r) &= \sum_{s=-\lfloor Nh \rfloor}^{\lfloor Nh \rfloor} w(t_s, t_r)\,\beta_i\left(t_r + \frac{s}{N}\right) \\
&= \sum_{s=-\lfloor Nh \rfloor}^{\lfloor Nh \rfloor} w(t_s, t_r)\left\{\beta_i(t_r) + \beta_i'(t_r)\frac{s}{N} + \frac{1}{2}\beta_i''(\tau_s)\frac{s^2}{N^2}\right\} \\
&= \beta_i(t_r) + \sum_{s=-\lfloor Nh \rfloor}^{\lfloor Nh \rfloor} w(t_s, t_r)\frac{1}{2}\beta_i''(\tau_s)\frac{s^2}{N^2},
\end{aligned}$$

where $|\tau_s - t_r| \leq \frac{s}{N}$.

Hence, by (RC1),

$$\left|\tilde{\beta}_i(t_r) - \beta_i(t_r)\right| \leq \frac{(Nh)^2}{N^2} C \sum_{s=-\lfloor Nh \rfloor}^{\lfloor Nh \rfloor} w(t_s, t_r) = Ch^2.$$

Therefore

$$\|\tilde{\beta} - \beta\|_2^2 \leq |\mathcal{A}_n|\,C^2 h^4, \tag{A.3}$$

for $h = h_n$ small enough.

If we choose $h_n \asymp \log(n)^{1/5} n^{-1/5} N^{-1/5}$, all terms in (A.1) are of the same order.

□

**Proof of Proposition 2.2**

We use the decomposition in (A.1). Since the variance in the smoothed case is of order $1/(Nh)$, see (A.2), we obtain

$$\|\hat{\tilde{\beta}} - \tilde{\beta}\|_2^2 \leq O_P((Nh)^{-1/2} a_n). \quad (A.4)$$

On the other hand, we have by (A.3)

$$\|\tilde{\beta} - \beta\|_2^2 \asymp |\mathcal{A}_n| h^4. \quad (A.5)$$

The optimal rate for the bandwidth minimizing the terms in (A.4) and (A.5) is

$$h_{opt} = N_n^{-1/9} |\mathcal{A}_n|^{-2/9} a_n^{2/9} \to \infty, \ n \to \infty$$

and we obtain using (A.1), (A.4) and (A.5)

$$\|\hat{\tilde{\beta}} - \beta\|_2^2 \leq O_P((Nh_{opt})^{-1/2} a_n). \quad (A.6)$$

Since

$$Nh_{opt} \asymp N_n^{8/9} |\mathcal{A}_n|^{-2/9} a_n^{2/9} \to \infty, \ n \to \infty$$

because $N_n \gg |\mathcal{A}_n|^{1/4} a_n^{-1/4}$, we see from (A.6) that a faster convergence rate $o_P(a_n)$ is achieved.

□

**Proof of Theorem 3.1**

As in the proof of Proposition 2.1, we have the model

$$\tilde{y} = X\tilde{\beta} + \tilde{\varepsilon}_N$$

for the smoothed response $\tilde{y}$. Multiplying both sides with $\sqrt{Nh}$ results in

$$\tilde{\tilde{y}} = \tilde{X}\tilde{\beta} + \tilde{\tilde{\varepsilon}}_N, \quad (A.7)$$

with $\tilde{\tilde{y}} = \sqrt{Nh}\tilde{y}$, $\tilde{X} = \sqrt{Nh}X$ and $\tilde{\tilde{\varepsilon}}_N = \sqrt{Nh}\tilde{\varepsilon}_N$.

Note that the variance of the error term $\tilde{\tilde{\varepsilon}}_N$ depends on $N$. As can be seen from (A.2), we have for $N \to \infty$

$$\text{Var}\left((\tilde{\tilde{\varepsilon}}_N)_i\right) \sim \sigma^2 \int K^2(x) dx.$$

Using the rescaled model (A.7), we can now adapt the proof of [24].

*L. Meier and P. Bühlmann/Smoothing $\ell_1$-penalized estimators* 613Let us first focus on the problem on the original scale. We re-parameterize the parameter vector $\beta$ as
$$\beta = \tilde{\beta} + \frac{u}{\sqrt{nNh}},$$
or $u = \sqrt{nNh}(\beta - \tilde{\beta}) \in \mathbb{R}^p$. The quantity of interest is $\hat{u} = \sqrt{nNh}(\hat{\tilde{\beta}} - \tilde{\beta})$, where
$$\hat{u} = \arg\min_u \psi_n(u),$$
with
$$\psi_n(u) = \left\| \tilde{y} - \sum_{j=1}^p x_j \left( \tilde{\beta}_j + \frac{u_j}{\sqrt{nNh}} \right) \right\|_2^2 + \lambda_n \sum_{j=1}^p \hat{w}_j \left| \tilde{\beta}_j + \frac{u_j}{\sqrt{nNh}} \right|.$$

By multiplying $\psi_n(u)$ with $Nh$, we can rewrite $\hat{u} = \arg\min_u \tilde{\psi}_n(u)$, where
$$\tilde{\psi}_n(u) = \left\| \tilde{\tilde{y}} - \sum_{j=1}^p \tilde{x}_j \left( \tilde{\beta}_j + \frac{u_j}{\sqrt{nNh}} \right) \right\|_2^2 + \tilde{\lambda}_n \sum_{j=1}^p \hat{w}_j \left| \tilde{\beta}_j + \frac{u_j}{\sqrt{nNh}} \right|,$$

and $\tilde{\lambda}_n = Nh\lambda_n$. Now we can follow the proof of [24]. With slight changes, because of the non-constant variance of the error-term, we arrive at
$$\sqrt{nNh}(\hat{\tilde{\beta}}_\mathcal{A} - \tilde{\beta}_\mathcal{A}) \xrightarrow{d} N(0, \sigma_*^2 (C_{\mathcal{A}\mathcal{A}})^{-1}),$$

where $\mathcal{A}$ is the active set of the *unsmoothed* parameter vector, i.e. the parameter vector at the current time-point. Finally, observe that
$$\sqrt{nNh}(\hat{\tilde{\beta}}_\mathcal{A} - \beta_\mathcal{A}) = \sqrt{nNh}(\hat{\tilde{\beta}}_\mathcal{A} - \tilde{\beta}_\mathcal{A}) + \sqrt{nNh}(\tilde{\beta}_\mathcal{A} - \beta_\mathcal{A}),$$

and that we get for the second term analogously as in (A.3), using $|\mathcal{A}_n| \leq p < \infty$,
$$nNh\|\tilde{\beta}_\mathcal{A} - \beta_\mathcal{A}\|_2^2 \leq CnNh^5$$

for some constant $C$. If we choose $h = o(n^{-1/5}N^{-1/5})$, the asymptotic normality part follows.

The proof of model selection consistency is analogous to [24]. □

## References

[1] BEER, M. A. AND TAVAZOIE, S. (2004). Predicting gene expression from sequence. *Cell 117*, 185–198.